\documentclass[11pt,letterpaper]{article}
\usepackage{fullpage}

\usepackage{times}
\usepackage{amsthm}
\usepackage{amsmath}
\usepackage{amssymb}
\newcommand{\comment}[1]{}

\newtheorem{theorem}{Theorem}

\newtheorem{cor}{Corollary}

\title{Values of $L$-Series of Hecke Eigenforms}
\author{Kamal Khuri-Makdisi, Winfried Kohnen and Wissam Raji}
\date{}
\begin{document}
\maketitle
\begin{abstract}
We determine a formula for the average values of $L$-series associated to eigenforms at complex values.
\end{abstract}

\section{Introduction}
Let $S_k$ denote the space of cusp forms of integer weight $k$ on the full modular group $SL_2(\mathbb{Z})$.  We define the $n$-th period of $f$ in $S_k$ by 

$$r_n(f)=\int_0^{\infty}f(it)t^ndt, \ \ \ 0\leq n\leq k-2.$$ It is well-known that $S_k$ has a rational structure given by the rationality of its periods \cite{KZ}.  Periods of cusp forms are actually the critical values of their corresponding $L$-series, and by definition also are the coefficients of the period polynomial associated to the modular form of degree $k-2$.  Petersson gives an average trace formula for the product of Fourier coefficients of cusp forms in terms of Kloosterman sums (see \cite{HI}).  In this paper, we give an analogue for Petersson's average formula where the Fourier coefficients are replaced by $L$-values of Hecke eigenforms at arbitrary values with certain restrictions.  Restricting to integers though, we obtain an average result of the explicit formulas proved by Kohnen and Zagier in \cite{KZ} for the periods of the kernel function for the special $L$-values.

Let us be more precise now. As mentioned above, the periods $r_n(f)$ give rise to a rational structure for the space of cusp forms via the Eichler-Shimura theory and those periods can be determined by taking the Petersson product of a certain kernel function $R_n$ with $f$ given by $r_n(f)=<f,R_n>$ for every $f\in S_k$.  In \cite{KZ}, the authors determine the periods of the function $R_n$ in terms of Bernoulli numbers and show certain symmetric properties of these periods. As a result, one can determine average values of L-series associated to Hecke eigenforms at integer values.  From what is mentioned, one can notice the importance of critical values of $L$-series at integer values. In this paper, we generalize the result from \cite{KZ} and determine a formula for the average values of the $L$-functions associated to Hecke eigenforms at complex values. The generalization is not immediate due to the branching problems that emerge and due to some complications in certain contour integrals.  It is worth mentioning that for general $s$ and $s'$ in $\mathbb{C}$, Theorem 1.2 of \cite{D} expresses a similar expression to (\ref{int2}) as the inner product of the first Poincare series $P_1$ with a product of two non-holomorphic Eisenstein series.

Let $k\in 2\mathbb{Z}$.  If $z\in \mathbb{C}-\{0\}$ and $w\in \mathbb{C}$, let $z^w=e^{w\log z}$ where $-\pi<arg z\leq \pi$.  If $f(z)=\sum_{n\geq 1}a(n)e^{2\pi inz}\in S_k$, we define the normalized L-series associated to $f$ by
\begin{equation}
L^*(f,s):=(2\pi)^{-s}\Gamma(s)L(f,s).
\end{equation} 
$L^*(f,s)$ has analytic continuation to $\mathbb{C}$ and satisfies the functional equation
\begin{equation}
L^*(f,k-s)=(-1)^{k/2}L^*(f,s).
\end{equation}
Let $\{f_{k,1},f_{k,2},...,f_{k,g_k}\}$ be the basis of normalized Hecke eigenforms of $S_k$.

We define the kernel function on $SL_2(\mathbb{Z})$ of integer weight $k$ as given in \cite{KZ} given by 
\begin{equation}\label{kernel}
R_{s,k}(z):=\gamma_{k}(s)\sum_{V} (cz+d)^{-k}\left(\frac{az+b}{cz+d}\right)^{-s}
\end{equation}
for $s\in \mathbb{C}$, $1<Re \ s< k-1$.  Here
the sum runs over all matrices $V=\left(%
\begin{array}{cc}
  a & b \\
  c & d \\
\end{array}%
\right)$ in $PSL_2(\mathbb{Z})$, and $$\gamma_k(s)=e^{\pi is/2}\Gamma(s)\Gamma(k-s).$$ 
In Lemma 1 in \cite{K1}, it was shown that under the usual Petersson inner product, one has for any $f\in S_k$
\begin{equation}
<f,R_{\bar{s},k}>=c_kL^*(f,s),
\end{equation}
where $$c_k=\frac{(-1)^{k/2}\pi (k-2)!}{2^{k-2}}.$$
By Lemma 1 in \cite{K1},

\begin{equation}\label{inner}
R_{s,k}=c_k\sum_{\nu=1}^{g_k}\frac{L^*(f_{k,\nu},s)}{<f_{k,\nu},f_{k,\nu}>}f_{k,\nu}.
\end{equation}

Taking the Mellin transform of both sides, we obtain for $s'\in \mathbb{C}$

\begin{equation}\label{int2}
L^*(R_{s,k},s')=\int_0^{\infty}R_{s,k}(it)t^{s'-1}dt=c_k\sum_{\nu=1}^{g_k}\frac{L^*(f_{k,\nu},s)L^*(f_{k,\nu},s')}{<f_{k,\nu},f_{k,\nu}>}.
\end{equation}





In what follows, we calculate the integral of the kernel function and deduce the value of the right hand side of (\ref{int2}).  The result in the following theorem is a generalization of part of the result of \cite{KZ}, due to the conditions on $s$ and $s'$.
\begin{theorem}\label{thm1}[The Main Theorem]
Let $s+s'\in 2\mathbb{Z}+1, Re \ s>Re \ s'+1, 1< s+ s'<k-1, 1< Re \ s < k-1,$ and $1<Re \ s'<k-1$.  Then 
\begin{equation}\label{mainint}
\begin{aligned}
\gamma_k(s)^{-1}\int_0^{\infty}R_{s,k}(it)t^{s'-1}dt&=  i^{-k}\frac{(2\pi i)^{k-s}}{(2\pi)^{k-s'}}\zeta(s-s'+1)\frac{\Gamma(k-s')}{\Gamma(k-s)}+i^{-k}\frac{(-2\pi i)^{s}}{(2\pi)^{k-s'}}\zeta(k-s'-s+1)\frac{\Gamma(k-s')}{\Gamma(s)}\\ &+2\pi  e^{\pi i (s'-1)/2}\frac{(-s'+1)_{k-1}}{(k-1)!} \sum_{\substack{a,b,c,d>0 \\ ad-bc=1}}a^{-s}c^{s-s'}d^{s'-k}\  {_2}F_1(s,-s'+k;k;1/ad)\\&+\zeta(s-s')\left(e^{-\pi is'/2}-e^{ 3\pi i s'/2}\right)\left(\frac{\Gamma(s')\Gamma(s-s')}{\Gamma(s)}\right),
\end{aligned}
\end{equation}
where ${_2}F_1$ is the hypergeometric series.
\end{theorem}

Note that $Re \ s>Re \ s'+1$ will be used in a crucial way in e.g. equation $(28)$ and therefore this condition cannot be replaced by the symmetric condition $Re (s-s')\not\in [-1,1]$. It is also worth noting that the condition that $s+s'$ is an odd integer will play a significant role in equation $(16)$ where we were able to recombine the integrals and get an integral over the real line.  However, that condition could be removed but one would then have to re-evaluate additional integrals and the strategy used to evaluate $S_{\epsilon}$ by considering equation (18) would need to be completely changed. We would have liked to remove the condition on $s$ and $s'$ since the expression in (6) is symmetric with respect to $s$ and $s'$ and it would be useful to be able to see this in (\ref{mainint}).  

\par We would like to point out that if we restrict the values of $s$ and $s'$ to integers, the third term on the right side of equation (\ref{mainint}) containing the hypergeometric series will become zero. In the last section, we show that this term emerges from evaluating an integral that was zero in the setting of \cite{KZ}. In our context, the factor $(-s'+1)_{k-1}$ becomes zero if $s'$ (and hence $s$) is an integer as in \cite{KZ}. Moreover, the remaining terms of (\ref{mainint}) will give the Bernoulli numbers in terms of the special values of the Riemann $\zeta$-function.

Using (\ref{int2}), we deduce the following corollary about the average values of the $L$-series.

\begin{cor}
For $s+s'\in 2\mathbb{Z}+1, 1< s+ s'<k-1, 1< Re \ s< k-1, 1<Re \ s'<k-1$  and $Re \ s> Re \ s'+1$, we have 
\begin{equation}
\begin{aligned}
\sum_{\nu=1}^{g_k}\frac{L^*(f_{k,\nu},s)L^*(f_{k,\nu},s')}{<f_{k,\nu},f_{k,\nu}>}&=\frac{\gamma_k(s)}{c_k}i^{-k}\frac{(2\pi i)^{k-s}}{(2\pi)^{k-s'}}\zeta(s-s'+1)\frac{\Gamma(k-s')}{\Gamma(k-s)}\\ &+\frac{\gamma_k(s)}{c_k}i^{-k}\frac{(-2\pi i)^{s}}{(2\pi)^{k-s'}}\zeta(k-s'-s+1)\frac{\Gamma(k-s')}{\Gamma(s)} \\&+\frac{2\pi \gamma_k(s)}{c_k} e^{\pi i (s'-1)/2}\frac{(-s'+1)_{k-1}}{(k-1)!} \sum_{\substack{a,b,c,d>0 \\ ad-bc=1}}a^{-s}c^{s-s'}d^{s'-k}\  {_2}F_1(s,-s'+k;k;1/ad)\\&+\frac{\gamma_k(s)}{c_k}\zeta(s-s')\left(e^{-\pi is'/2}-e^{ 3\pi i s'/2}\right)\left(\frac{\Gamma(s')\Gamma(s-s')}{\Gamma(s)}\right).\end{aligned}
\end{equation}
\end{cor}
For $k=8,10,14$ one has $S_k= \{0\}$ and we deduce a relation between $\zeta(s)$ and the hypergeometric function ${_2}F_1$. Note that for $k=4$ and $k=6$, there are no values for $s$ and $s'$ that satisfy the conditions.

\begin{cor} For $k=8,10, 14$ and for $s+s'\in 2\mathbb{Z}+1, 1< s+ s'<k-1, 1< Re \ s< k-1, 1<Re \ s'<k-1$  and $Re \ s>Re \ s' +1$, we have 
\begin{equation}
\begin{aligned}
&i^{-k}\frac{(2\pi i)^{k-s}}{(2\pi)^{k-s'}}\frac{\Gamma(k-s')}{\Gamma(k-s)}\zeta(s-s'+1)+i^{-k}\frac{(-2\pi i)^{s}}{(2\pi)^{k-s'}}\frac{\Gamma(k-s')}{\Gamma(s)}\zeta(k-s'-s+1)\\&+\zeta(s-s')\left(e^{-\pi is'/2}-e^{ 3\pi i s'/2}\right)\left(\frac{\Gamma(s')\Gamma(s-s')}{\Gamma(s)}\right)\\&=-2\pi  e^{\pi i (s'-1)/2}\frac{(-s'+1)_{k-1}}{(k-1)!} \sum_{\substack{a,b,c,d>0 \\ ad-bc=1}}a^{-s}c^{s-s'}d^{s'-k}\  {_2}F_1(s,-s'+k;k;1/ad).
\end{aligned}
\end{equation}
\end{cor}
\section{Proof of the Main Theorem}
\begin{proof}
We follow a strategy similar to \cite{KZ}, by considering three cases pertaining to the matrices $V$ in the definition of the kernel function separately. In other words, we divide the sum of $R_{s,k}$ defined in (\ref{kernel}) inside the integral on the left side of equation (\ref{mainint}) into three parts according to whether $bd=0$, $bd>0$ and $bd<0$.  In what follows, we treat the case $bd=0$ separately, and then combine the cases $bd>0$ and $bd<0$ into one integral. We then use a contour integral evaluation to deal with convergence issues. We now start with the first case when $bd=0$.  

For $bd=0$, we have elements $V$ of the form $\left(%
\begin{array}{cc}
  1 & 0 \\
  n & 1 \\
\end{array}%
\right)$ and  $\left(%
\begin{array}{cc}
  n & -1 \\
  1 & 0 \\
\end{array}%
\right)$. 
  
Thus we have a contribution to $\gamma_k(s)^{-1}\int_0^{\infty}R_{s,k}(it)t^{s'-1}dt$ of the form
\begin{equation}
\begin{aligned}
A&:=\int_0^{\infty}\sum_{bd=0} (cit+d)^{-k}\left(\frac{ait+b}{cit+d}\right)^{-s}t^{s'-1}dt\\&= \int_0^{\infty}\sum_{n\in \mathbb{Z}}(n-i/t)^{s-k}(it)^{-k}t^{s'-1}dt +\int_0^{\infty}\sum_{n\in \mathbb{Z}}(n+i/t)^{-s}(it)^{-k}t^{s'-1}dt\\&=i^{-k}\int_0^{\infty}\sum_{n\in \mathbb{Z}}(n-i/t)^{s-k}t^{s'-k-1}dt +i^{-k}\int_0^{\infty}\sum_{n\in \mathbb{Z}}(n+i/t)^{-s}t^{s'-k-1}dt.
\end{aligned}
\end{equation}
We now replace $t$ by $1/t$ to obtain
\begin{equation}
A=i^{-k}\int_0^{\infty}\sum_{n\in \mathbb{Z}}(n-it)^{s-k}t^{k-s'-1}dt +i^{-k}\int_0^{\infty}\sum_{n\in \mathbb{Z}}(n+it)^{-s}t^{k-s'-1}dt.
\end{equation}
The Lipschitz summation formula tells us that for $s\in \mathbb{C}$, with $ Re \ s >0$,
$$\sum_{n\in \mathbb{Z}}(z+n)^{-s-1}=\frac{(-2\pi i)^{s+1}}{\Gamma(s+1)} \sum_{n\geq 1}n^{s}e^{2\pi inz}, \ \ \  \Im z>0.$$ Letting $z=it$ in the Lipschitz summation formula for real $s$ and then using the complex conjugate, we can analytically continue this formula for complex $s$.  We get
\begin{equation}
\begin{aligned}
A&= i^{-k}\frac{(2\pi i)^{k-s}}{\Gamma(k-s)}\int_0^{\infty}\sum_{n\geq 1}n^{k-s-1}e^{-2\pi nt}t^{k-s'-1}dt +i^{-k}\frac{(-2\pi i)^{s}}{\Gamma(s)}\int_0^{\infty}\sum_{n\geq 1}n^{s-1}e^{-2\pi nt}t^{k-s'-1}dt\\&= i^{-k}\frac{(2\pi i)^{k-s}}{(2\pi)^{k-s'}}\zeta(s-s'+1)\frac{\Gamma(k-s')}{\Gamma(k-s)}+i^{-k}\frac{(-2\pi i)^{s}}{(2\pi)^{k-s'}}\zeta(k-s'-s+1)\frac{\Gamma(k-s')}{\Gamma(s)}
\end{aligned}
\end{equation}
since $s+s'<k-2$ and $Re s> Re s'$ are satisfied.

As for the other elements where $bd<0$ and $bd>0$, we set
$$B:=\int_0^{\infty}\sum_{bd<0} (cit+d)^{-k}\left(\frac{ait+b}{cit+d}\right)^{-s}t^{s'-1}dt$$ and $$B':=\int_0^{\infty}\sum_{bd>0} (cit+d)^{-k}\left(\frac{ait+b}{cit+d}\right)^{-s}t^{s'-1}dt.$$
Replace $t$ by $-t$ to get
\begin{equation}
B=\int_{-\infty}^0\sum_{bd<0} (-cit+d)^{-k}\left(\frac{-ait+b}{-cit+d}\right)^{-s}(-t)^{s'-1}dt.
\end{equation} 
Now replace $(a,d)$ by $(-a,-d)$, so we have
\begin{equation}
B=\int_{-\infty}^0\sum_{bd>0} (-cit-d)^{-k}\left(\frac{ait+b}{-cit-d}\right)^{-s}(-t)^{s'-1}dt.
\end{equation}
As a result, we have 
\begin{equation}
B=e^{-\pi i(s+s'-1)}\int_{-\infty}^0\sum_{bd>0} (cit+d)^{-k}\left(\frac{ait+b}{cit+d}\right)^{-s} t^{s'-1}dt=\int_{-\infty}^0\sum_{bd>0} (cit+d)^{-k}\left(\frac{ait+b}{cit+d}\right)^{-s}t^{s'-1}dt
\end{equation}
since $s+s'\in 2\mathbb{Z}+1$. 
We combine the cases $bd>0$ and $bd<0$ to get the following sum 
\begin{equation}
B+B'=\int_0^{\infty}\sum_{bd>0} (cit+d)^{-k}\left(\frac{ait+b}{cit+d}\right)^{-s}t^{s'-1}dt+\int_{-\infty}^0\sum_{bd>0} (cit+d)^{-k}\left(\frac{ait+b}{cit+d}\right)^{-s}t^{s'-1}dt.
\end{equation}
To be able to interchange summation and integration, we note that the above series converge uniformly for $t$ in a compact subset of $\mathbb{R}-\{0\}$.  Let us therefore define the following integrals:
$$S_{\epsilon}:=\sum_{bd>0}\int_{\epsilon}^{\frac{1}{\epsilon}} (cit+d)^{-k}\left(\frac{ait+b}{cit+d}\right)^{-s}t^{s'-1}dt+ \sum_{bd>0}\int_{-\frac{1}{\epsilon}}^{-\epsilon}(cit+d)^{-k}\left(\frac{ait+b}{cit+d}\right)^{-s}t^{s'-1}dt.$$
Notice that the interchange of sums and integrals is justified in $S_{\epsilon}$, and hence $\lim_{\epsilon \rightarrow 0}S_{\epsilon}=B+B'.$   Once we calculate $\lim_{\epsilon \rightarrow 0}S_{\epsilon}$, we combine this with our value of $A$, and we obtain

$$\gamma_k(s)^{-1}\int_0^{\infty}R_{s,k}(it)t^{s'-1}dt=A+B+B'.$$

\subsection{The Value of $\lim_{\epsilon \rightarrow 0}S_{\epsilon}$}
In what follows, we calculate $\lim_{\epsilon \rightarrow 0}S_{\epsilon}$.
We have

\begin{equation}\label{intsum}
\sum_{bd>0}\left(\int_{\epsilon}^{1/\epsilon}+\int_{-1/\epsilon}^{-\epsilon}\right)=\sum_{bd>0}\left(
\int_{-\infty}^{\infty}-\int_{-\epsilon}^{\epsilon}
-\int_{-\infty}^{-1/\epsilon}-\int_{1/\epsilon}^{\infty}\right)=\sum_{bd>0}\left(\int_{-\infty}^{\infty}-\int_{-\epsilon}^{\epsilon}\right)+S_{\epsilon}^0,
\end{equation}
where $$S_{\epsilon}^0:=\sum_{bd>0}\left(-\int_{-\infty}^{-1/\epsilon}-\int_{1/\epsilon}^{\infty}\right).$$
We now try to simplify $S_{\epsilon}^0.$ Substitute $t$ by $1/t$ in $S_{\epsilon}^0$.  This gives

\begin{equation}
\begin{aligned}
S_{\epsilon}^0 &=\sum_{bd>0}\left(\int_{\epsilon}^0 (ci+dt)^{-k}\left(\frac{ai+bt}{ci+dt}\right)^{-s}t^{k-s'-1}dt+\int_0^{-\epsilon} (ci+dt)^{-k}\left(\frac{ai+bt}{ci+dt}\right)^{-s}t^{k-s'-1}dt\right) \\ &=-\sum_{bd>0}\int_{-\epsilon}^{\epsilon}(ci+dt)^{-k}\left(\frac{ai+bt}{ci+dt}\right)^{-s}t^{k-s'-1}dt.
\end{aligned}
\end{equation}

We now replace $\left(%
\begin{array}{cc}
  a & b \\
  c & d \\
\end{array}%
\right)$ by  $\left(%
\begin{array}{cc}
  b & -a \\
  d &  -c\\
\end{array}%
\right)$, and get
\begin{equation}
\begin{aligned}
S_{\epsilon}^0 &=-i^{-k}\sum_{ac>0}\int_{-\epsilon}^{\epsilon}(cit+d)^{-k}\left(\frac{ait+b}{cit+d}\right)^{-s}t^{k-s'-1}dt.
\end{aligned}
\end{equation}
 Combining with (\ref{intsum}), we obtain
\begin{equation}\label{sum3}
\begin{aligned}
S_{\epsilon}&=\sum_{bd>0}\int_{-\infty}^{\infty} (cit+d)^{-k}\left(\frac{ait+b}{cit+d}\right)^{-s}t^{s'-1}dt-\sum_{bd>0}\int_{-\epsilon}^{\epsilon}(cit+d)^{-k}\left(\frac{ait+b}{cit+d}\right)^{-s}t^{s'-1}dt\\&-i^{-k}\sum_{ac>0}\int_{-\epsilon}^{\epsilon}(cit+d)^{-k}\left(\frac{ait+b}{cit+d}\right)^{-s}t^{k-s'-1}dt.
\end{aligned}
\end{equation}

We have $ac.bd=ad.bc=(1+bc)bc\geq 0$, which implies that $bd\geq 0$ in the third sum in (\ref{sum3}). Hence
\begin{equation}\label{sum4}
\begin{aligned}
S_{\epsilon}&=\sum_{bd>0}\int_{-\infty}^{\infty} (cit+d)^{-k}\left(\frac{ait+b}{cit+d}\right)^{-s}t^{s'-1}dt-i^{-k}\sum_{ac>0, bd=0}\int_{-\epsilon}^{\epsilon}(cit+d)^{-k}\left(\frac{ait+b}{cit+d}\right)^{-s}t^{k-s'-1}dt\\&+i^{-k}\sum_{ac=0, bd>0}\int_{-\epsilon}^{\epsilon}(cit+d)^{-k}\left(\frac{ait+b}{cit+d}\right)^{-s}t^{k-s'-1}dt\\ &-\sum_{bd>0}\left(\int_{-\epsilon}^{\epsilon}(cit+d)^{-k}\left(\frac{ait+b}{cit+d}\right)^{-s}t^{s'-1}dt+i^{-k}\int_{-\epsilon}^{\epsilon}(cit+d)^{-k}\left(\frac{ait+b}{cit+d}\right)^{-s}t^{k-s'-1}dt\right). 
\end{aligned}
\end{equation}
Now to determine $S_{\epsilon}$, we divide our integral as follows.  Put $$S_{\epsilon}=S+S_{\epsilon}'+S_{\epsilon}''+S_{\epsilon}''',$$ where
\begin{equation}\label{S}
S:=\sum_{bd>0}\int_{-\infty}^{\infty} (cit+d)^{-k}\left(\frac{ait+b}{cit+d}\right)^{-s}t^{s'-1}dt,
\end{equation}
\begin{equation}\label{S'}
S_{\epsilon}':=-i^{-k}\sum_{ac>0, bd=0}\int_{-\epsilon}^{\epsilon}(cit+d)^{-k}\left(\frac{ait+b}{cit+d}\right)^{-s}t^{k-s'-1}dt,
\end{equation}
\begin{equation}\label{S''}
S_{\epsilon}'':=i^{-k} \sum_{ac=0, bd>0}\int_{-\epsilon}^{\epsilon}(cit+d)^{-k}\left(\frac{ait+b}{cit+d}\right)^{-s}t^{k-s'-1}dt
\end{equation}
and
\begin{equation}\label{S'''}
 S_{\epsilon}''':=-\sum_{bd>0}\left(\int_{-\epsilon}^{\epsilon}(cit+d)^{-k}\left(\frac{ait+b}{cit+d}\right)^{-s}t^{s'-1}dt+i^{-k}\int_{-\epsilon}^{\epsilon}(cit+d)^{-k}\left(\frac{ait+b}{cit+d}\right)^{-s}t^{k-s'-1}dt\right). 
 \end{equation}
 The evaluation of $S$ is given in the last section of this paper. We show there that  
 
\begin{equation}\label{valS}
\begin{aligned}
S&=2\pi  e^{\pi i (s'-1)/2}\frac{(-s'+1)_{k-1}}{(k-1)!} \sum_{\substack{bd>0 \\ ad-bc=1\\ac\neq 0}}a^{-s}c^{s-s'}d^{s'-k}\  {_2}F_1(s,-s'+k;k;1/ad)\\&+i\zeta(s-s')\left(e^{3\pi i(s'-1)/2}-e^{-\pi i(s'-1)/2}\right)\left(\frac{\Gamma(s')\Gamma(s-s')}{\Gamma(s)}\right).
\end{aligned}
\end{equation}
 
 We now evaluate the remaining limits $\lim_{\epsilon \rightarrow 0} S'_{\epsilon}$, $\lim_{\epsilon \rightarrow 0 }S''_{\epsilon}$ and $\lim_{\epsilon \rightarrow 0} S'''_{\epsilon}$.
We start with $\lim_{\epsilon \rightarrow 0} S'_{\epsilon}$ defined in (\ref{S'}). Separate the sum into two sums where $b=0$ and where $d=0$. Thus 
\begin{equation}
S'_{\epsilon}=-i^{-k}\sum_{n>0}\int_{-\epsilon}^{\epsilon}(nit+1)^{-k}\left(\frac{it}{nit+1}\right)^{-s}t^{k-s'-1}dt-i^{-k}\sum_{n>0}\int_{-\epsilon}^{\epsilon}(it)^{-k}\left(\frac{nit-1}{it}\right)^{-s}t^{k-s'-1}dt.
\end{equation}
Substituting $t$ by $\epsilon t$ and $t$ by $-\epsilon t$ respectively in the above sums, we get
\begin{equation}\label{29}
\begin{aligned}
S'_{\epsilon}&=-i^{-k}\epsilon^{k-s'-s-1}\left(\epsilon \sum_{n=1}^{\infty}\int_{-1}^1(ni\epsilon t+1)^{-k}\left(\frac{it}{ni\epsilon t+1}\right)^{-s}t^{k-s'-1}dt\right)\\&-e^{\pi i(s-s'+1)}\epsilon^{s-s'-1}\left(\epsilon \sum_{n=1}^{\infty}\int_{-1}^1(it)^{-k}\left(\frac{ni\epsilon t-1}{it}\right)^{-s}t^{k-s'-1}dt \right).
\end{aligned}
\end{equation} 
We now take the limit as $\epsilon \rightarrow 0$ of the two summands.  Observe that the two expressions in brackets are Riemann sums for the integrals 
$$\int_{0}^{\infty}\int_{-1}^1(ixt+1)^{-k}\left(\frac{it}{ixt+1}\right)^{-s}t^{k-s'-1}dtdx\ \ \mbox{and} \ \ \ \int_{0}^{\infty}\int_{-1}^1(it)^{-k}\left(\frac{ixt-1}{it}\right)^{-s}t^{k-s'-1}dt dx.$$

Since $s+s'<k-1$ and $Re \ s- Re \ s'> 1$, the powers of $\epsilon$ to the left of the parentheses allow us to conclude that
\begin{equation}\label{valS'}
\lim_{\epsilon \rightarrow 0}S'_{\epsilon}=0
\end{equation}

As for $S_{\epsilon}''$ as defined in (\ref{S''}), we follow a similar evaluation as for $\lim_{\epsilon \rightarrow 0}S_{\epsilon}'$ as above, to show that the $\lim_{\epsilon\rightarrow 0} S_{\epsilon}''=0$.  By directly estimating the integrals that emerge from the cases where $a=0$ and $c=0$, one can easily check that
\begin{equation}\label{valS''}
\lim_{\epsilon \rightarrow 0}S''_{\epsilon}=0.
\end{equation}
For example, in the case $a=0$, replace $t$ by $\epsilon t$, so the summand is
\begin{equation*}
\epsilon^{k-s'}\sum_{n>0}\int_{-1}^{1}\left(\frac{1}{-i\epsilon t+n}\right)^{-s}(-i\epsilon t+n)^{-k}t^{k-s'-1}dt.
\end{equation*}
Note that 
$$\big|\sum_{n>0}\int_{-1}^{1}\left(\frac{1}{-i\epsilon t+n}\right)^{-s}(-i\epsilon t+n)^{-k}t^{k-s'-1}dt\big|\leq \sum_{n> 0}n^{Re \ s-k}\int_{-1}^1t^{k-Re \ s'-1}dt.$$
The sum above converges because $Re\ s-k<-1$ and the integral converges because $k-Re\ s'-1>-1$. 
 A similar computation covers the case $c=0$.  This completes the proof of (30).\\ 
What is left is to determine the $\lim_{\epsilon\rightarrow 0} S_{\epsilon}'''$.  Recall that our group is $PSL_2(\mathbb{Z})$ and that the sums below run over all $bd>0$. Without loss of generality, we can take $b>0$ and $d>0$.  Choose $a_0$ and $c_0$ such that $a_0d-bc_0=1$. All the other elements with the same $(b,d)$ are given as $\left(%
\begin{array}{cc}
  a_0+bn & b \\
  c_0+dn & d \\
\end{array}%
\right).$ Now taking $l=\frac{c_0}{d}+n=\frac{c}{d}$, one notices that the first summand in $S_{\epsilon}'''$ can be written as
\begin{equation}
\begin{aligned}
&-\sum_{bd>0}\int_{-\epsilon}^{\epsilon}(cit+d)^{-k}\left(\frac{ait+b}{cit+d}\right)^{-s}t^{s'-1}dt \\&=-\sum_{b,d>0, (b,d)=1}\frac{1}{d^{k-s}b^s}\int_{-\epsilon}^{\epsilon}\sum_{l\in \mathbb{Z}+c_0/d}(ilt+1)^{-k}\left(\frac{1+it/bd+ilt}{ilt+1}\right)^{-s}t^{s'-1}dt.
\end{aligned}
\end{equation}
We replace $t$ by $\epsilon t$ and get
\begin{equation}
\begin{aligned}
&-\sum_{b,d>0}\int_{-\epsilon}^{\epsilon}(cit+d)^{-k}\left(\frac{ait+b}{cit+d}\right)^{-s}t^{s'-1}dt\\&=-\epsilon^{s'-1}\sum_{b,d>0, (b,d)=1}\frac{1}{d^{k-s}b^s}\left(\epsilon \sum_{l\in \mathbb{Z}+c_0/d} \int_{-1}^{1}(il\epsilon t+1)^{-k}\left(\frac{1+i\epsilon t/bd+il\epsilon t}{il\epsilon t+1}\right)^{-s}t^{s'-1}dt\right).
\end{aligned}
\end{equation}
Similarly, the expression in brackets above gives a Riemann sum for an integral. So we obtain $\epsilon^{s'-1}$ times a quantity whose limit is 
\begin{equation}
-\sum_{b,d>0, (b,d)=1}d^{s-k}b^{-s}\int_{-\infty}^{\infty}\int_{-1}^{1}\frac{dt}{(1+ixt)^k}dx=\frac{4\pi \zeta(k-s)\zeta(s)}{(k-1)\zeta(k)}. 
\end{equation}
Note that $Re \ s'>1$ and thus the expression vanishes as $\epsilon \rightarrow 0$.

We similarly deal with the second sum of $S_{\epsilon}'''$ as defined in (\ref{S'''}).  Here the limit vanishes as $\epsilon \rightarrow 0$ because $Re \ s' <k-1$. As a result, 
\begin{equation}\label{valS'''}
\lim_{\epsilon \rightarrow 0}S_{\epsilon}'''=0.
\end{equation}
Adding eqs (\ref{valS}),  (\ref{valS'}), (\ref{valS''}) and  (\ref{valS'''}), we get
\begin{equation}
\begin{aligned}
B=\lim_{\epsilon \rightarrow 0}S_{\epsilon}&=2\pi  e^{\pi i (s'-1)/2}\frac{(-s'+1)_{k-1}}{(k-1)!} \sum_{\substack{a,b,c,d>0 \\ ad-bc=1}}a^{-s}c^{s-s'}d^{s'-k}\  {_2}F_1(s,-s'+k;k;1/ad)\\&+i\zeta(s-s')\left(e^{3\pi i(s'-1)/2}-e^{-\pi i(s'-1)/2}\right)\left(\frac{\Gamma(s')\Gamma(s-s')}{\Gamma(s)}\right).
\end{aligned}
\end{equation}
Finally, 

\begin{equation}
\begin{aligned}
\gamma_k(s)^{-1}\int_0^{\infty}R_{s,k}(it)t^{s'-1}dt&=  i^{-k}\frac{(2\pi i)^{k-s}}{(2\pi)^{k-s'}}\zeta(s-s'+1)\frac{\Gamma(k-s')}{\Gamma(k-s)}+i^{-k}\frac{(-2\pi i)^{s}}{(2\pi)^{k-s'}}\zeta(k-s'-s+1)\frac{\Gamma(k-s')}{\Gamma(s)}\\ &+2\pi  e^{\pi i (s'-1)/2}\frac{(-s'+1)_{k-1}}{(k-1)!} \sum_{\substack{a,b,c,d>0 \\ ad-bc=1}}a^{-s}c^{s-s'}d^{s'-k}\  {_2}F_1(s,-s'+k;k;1/ad)\\&+i\zeta(s-s')\left(e^{3\pi i(s'-1)/2}-e^{-\pi i(s'-1)/2}\right)\left(\frac{\Gamma(s')\Gamma(s-s')}{\Gamma(s)}\right).
\end{aligned}
\end{equation}
\end{proof}
\section{The evaluation of $S$}
In this section, we evaluate $S$ defined in subsection 2.1.   We start with the terms when $ac\neq 0$.  Note that since $bd>0$, we can assume that $b$ and $d$ are both positive.  As a result, by using the same argument between equations (\ref{sum3}) and (\ref{sum4}), we have that $ac>0$.  If $a$ and $c$ are both negative while $b$ and $d$ are both positive, we replace $\int_{-\infty}^{\infty}$ by the standard contour in the upper half plane that joins some small $\delta>0$ to some large $R>0$ on the real positive axis, the semi-circle of radius $R$ in the upper half plane and the segment above the negative real axis that joins $-R$ to $-\delta$ to skip the cut and then joins $-\delta$ to $\delta$ in the upper half plane. The singularities of 
$$\int_{-\infty}^{\infty} (cit+d)^{-k}\left(\frac{ait+b}{cit+d}\right)^{-s}t^{s'-1}dt$$ are at $t=ib/a$ and $t=id/c$ and in this case, they both lie in the lower half plane.  Moreover, the integrand on the semi-circle is $O(R^{Re \ s' -k})$.  As a result, the original integral is $0$, because $Re\ s'<k-1$.  \par On the other hand, if $a,b,c,d$ are all positive, then we have the poles $t=ib/a$ and $t=id/c$ in the upper half plane. Note that the integral is single-valued in the upper half plane minus the cut joining $ib/a$ and $id/c$. In this case, we take the same contour as the previous case when $a,c < 0$, which can be deformed to the contour around the cut joining $t=ib/a$ and $t=id/c$.  We call $C$ the contour with positive orientation around the segment cut joining $t=ib/a$ and $t=id/c$.  As a result, we have

$$\int_{-\infty}^{\infty} (cit+d)^{-k}\left(\frac{ait+b}{cit+d}\right)^{-s}t^{s'-1}dt=\int_C(cit+d)^{-k}\left(\frac{ait+b}{cit+d}\right)^{-s}t^{s'-1}dt.$$
In order to evaluate the integral at $C$, we make a change of variables by taking $x=a(cit+d)$, so we have

$$\int_C(cit+d)^{-k}\left(\frac{ait+b}{cit+d}\right)^{-s}t^{s'-1}dt=\frac{1}{i}a^{k-s-1}c^{s-s'}d^{s'-1}e^{\pi i (s'-1)/2}\int_{C'}x^{-k}\left(\frac{x-1}{x}\right)^{-s}\left(1-\frac{x}{ad}\right)^{s'-1}dx$$ where $C'$ is a contour going counterclockwise around the segment $[0,1]$ that keeps close to the segment (eg., a ``dumbbell'' contour).  Note that $ad=1+bc\geq 2$, so we can expand the binomial series (which converges for $|x|<ad$) and obtain
\begin{equation}
\begin{aligned}
&\frac{1}{i}a^{k-s-1}c^{s-s'}d^{s'-1}e^{\pi i (s'-1)/2}\int_{C'}x^{-k}\left(\frac{x-1}{x}\right)^{-s}\left(1-\frac{x}{ad}\right)^{s'-1}dx\\&=\frac{1}{i}a^{k-s-1}c^{s-s'}d^{s'-1}e^{\pi i (s'-1)/2}\sum_{m\geq 0}\int_{C'}x^{m-k}\left(\frac{x-1}{x}\right)^{-s}\frac{(-s'+1)_m}{a^md^mm!}dx\
\end{aligned}
\end{equation}
where $(.)_m$ is the rising Pochhammer symbol.
Now that we have expanded the series, we can deform $C'$ to come close to $\infty$. We use the residue with respect to the pole at $\infty$, which we rewrite using the substitution $w=1/x$ to get 
\begin{equation}
\int_{C'}x^{m-k}\left(\frac{x-1}{x}\right)^{-s}dx=2\pi iRes_{w=0}\left( w^{-m+k-2}(1-w)^{-s}dw\right)=\frac{(s)_l}{l!},
\end{equation}
where $l=m-k+1$. 
As a result, we get
\begin{equation}
\begin{aligned}
&\frac{1}{i}a^{k-s-1}c^{s-s'}d^{s'-1}e^{\pi i (s'-1)/2}\sum_{m\geq 0}\int_{C'}x^{m-k}\left(\frac{1-x}{x}\right)^{-s}\frac{(-s'+1)_m}{a^md^mm!}dx
\\ &=2\pi a^{k-s-1}c^{s-s'}d^{s'-1}e^{\pi i (s'-1)/2}\sum_{l \geq 0}\frac{(s)_l(-s'+1)_{k-1+l}}{l! (k-1+l)!}\left(\frac{1}{ad}\right)^{k-1+l} \\&=2\pi  a^{-s}c^{s-s'}d^{s'-k}e^{\pi i (s'-1)/2}\frac{(-s'+1)_{k-1}}{(k-1)!} \  {_2}F_1(s,-s'+k;k;1/ad).
\end{aligned}
\end{equation}
Finally, if $ac=0$, let us start with the case $a=0$.  Recall that $bd>0$, hence we have to evaluate 
\begin{equation}
\sum_{n>0}\int_{-\infty}^{\infty}(-it+n)^{s-k}t^{s'-1}dt.
\end{equation}
Notice that one can consider the standard semi-circle contour in the upper half plane to evaluate this integral.  The poles of $(-it+n)^{-k+s}$ are at $t=-in$ and thus we have 
\begin{equation}
\sum_{n>0}\int_{-\infty}^{\infty}(-it+n)^{s-k}t^{s'-1}dt=0.
\end{equation}

On the other hand, if $c=0$, we take $t=n\tau$, and we have
\begin{equation}
\sum_{n>0}\int_{-\infty}^{\infty}(it+n)^{-s}t^{s'-1}dt=\zeta(s-s')\int_{-\infty}^{\infty}(i\tau+1)^{-s}\tau^{s'-1}d\tau.
\end{equation}

We now divide the integral and substitute  $\tau=t$ and $\tau=-v$ respectively to obtain

\begin{equation}\label{twoints}
\begin{aligned}
\zeta(s-s')\int_{-\infty}^{\infty}(i\tau+1)^{-s}\tau^{s'-1}d\tau&=\zeta(s-s')\left(\int_0^{\infty}(it+1)^{-s}t^{s'-1}dt+\int_{0}^{\infty}(-iv+1)^{-s}(-v)^{s'-1}dv\right)\\&=\zeta(s-s')\left(\int_0^{\infty}(it+1)^{-s}t^{s'-1}dt+e^{\pi i(s'-1)}\int_{0}^{\infty}(-iv+1)^{-s}v^{s'-1}dv\right).
\end{aligned}
\end{equation}
We now evaluate each of the above integrals separately. We start with $\int_{0}^{\infty}(-iv+1)^{-s}v^{s'-1}dv$.  Making the change of variables $v=iw$, we get
\begin{equation}
\int_{0}^{\infty}(-iv+1)^{-s}v^{s'-1}dv=ie^{i\pi s/2}\int_{0}^{-i\infty}(w+1)^{-s}w^{s'-1}dw.
\end{equation}

One can easily show by estimating the integral over the contour in the fourth quadrant consisting of the quarter circle connecting $-i\infty$ to $+\infty$ that 
$$ie^{i\pi s/2}\int_{0}^{-i\infty}(w+1)^{-s}w^{s'-1}dw=ie^{i\pi s/2}\int_{0}^{\infty}(w+1)^{-s}w^{s'-1}dw.$$
Recall that for $Re \ p>0$ and $Re \ q>0$, we have the beta-integral given by

$$B(p,q)=\int_0^1x^{p-1}(1-x)^{q-1}dx=\frac{\Gamma(p)\Gamma(q)}{\Gamma(p+q)},$$ where $\Gamma(x)$ is the usual Gamma function.
A standard substitution (eg. (8.7-4) of \cite{H}) gives
\begin{equation}\label{oneint}
ie^{i\pi s/2}\int_{0}^{\infty}(w+1)^{-s}w^{s'-1}dw=ie^{i\pi (s'-1)/2}\frac{\Gamma(s')\Gamma(s-s')}{\Gamma(s)}.
\end{equation}
The second integral $\int_0^{\infty}(it+1)^{-s}t^{s'-1}dt$ can be evaluated using the same method, but now taking $w=it$ and the quarter circle contour in the first quadrant.  So we get
\begin{equation}\label{twoint}
\int_0^{\infty}(it+1)^{-s}t^{s'-1}dt=-ie^{-\pi i(s'-1)/2}\int_{0}^{\infty}(v+1)^{-s}v^{s'-1}dv=-ie^{-\pi i(s'-1)/2}\frac{\Gamma(s')\Gamma(s-s')}{\Gamma(s)}.
\end{equation}
Substituting (\ref{oneint}) and (\ref{twoint}) in $(\ref{twoints})$, we have
\begin{equation}
\zeta(s-s')\int_{-\infty}^{\infty}(i\tau+1)^{-s}\tau^{s'-1}d\tau=i\zeta(s-s')\left(e^{3\pi i(s'-1)/2}-e^{-\pi i(s'-1)/2}\right)\left(\frac{\Gamma(s')\Gamma(s-s')}{\Gamma(s)}\right).
\end{equation}

As a result, we have 
\begin{equation}
\begin{aligned}
S&=2\pi  e^{\pi i (s'-1)/2}\frac{(-s'+1)_{k-1}}{(k-1)!} \sum_{\substack{bd>0 \\ ad-bc=1; ac\neq 0}}a^{-s}c^{s-s'}d^{s'-k}\  {_2}F_1(s,-s'+k;k;1/ad)\\ &+i\zeta(s-s')\left(e^{3\pi i(s'-1)/2}-e^{-\pi i(s'-1)/2}\right)\left(\frac{\Gamma(s')\Gamma(s-s')}{\Gamma(s)}\right) .
\end{aligned}
\end{equation}

\newpage
\bibliographystyle{plain}

\begin{thebibliography}{WWW99}
\bibitem{D} N. Diamantis  and C. O'Sullivan , \textit{Kernels of L-functions of cusp forms}, Mathematische Annalen. 346(4) (2010) 897-929.
\bibitem{H} P. Henrici, \textit{Applied and Computational Complex Analysis}, volume 2, Wiley Classics Library, 1974.
\bibitem{HI} H. Iwaniec, \textit{Spectral methods of automorphic forms}, volume 53, American Mathematical
Society Providence, 2002.
\bibitem{K} W. Kohnen, \textit{Modular forms of half-integral weight on $\Gamma_0(4)$ }, Math Ann. 248 (1980), 249-266.
\bibitem{K1} W. Kohnen, \textit{Non-vanishing of Hecke L-functions},  J. of Number theory 67 (1997), 182-189.
\bibitem{KZ}
W.\ Kohnen and D.\ Zagier,
\newblock \textit{Modular forms with rational periods.}
\newblock In: Modular Forms (Durham, 1983), 1984, 197--249. 
\bibitem{Peter} H. Petersson, \textit{Uber die Entwicklungskoeffizienten der automorphen Formen}, Acta Math. 58,
169-215 (1932).
\bibitem{S} G. Shimura, \textit{On modular forms of half integral weight}, Annals of Math. 97 (1973), 440-481.

\end{thebibliography}

\small{Kamal Khuri-Makdisi, American University of Beirut, Beirut, Lebanon\\ \small \textit{Email address}: \textbf{kmakdisi@aub.edu.lb}}\\
\small{ Winfried Kohnen, University of Heidelberg, Im Neuenheimer Feld 205, 69120 Heidelberg, Germany}\\  \small \textit{E-mail address}: \textbf{winfried@mathi.uni-heidelberg.de}\\
{Wissam Raji, American University of Beirut, Beirut, Lebanon}\\  \small \textit{E-mail address}: \textbf{wr07@aub.edu.lb}

\end{document}